\documentclass[11pt]{article}   
\usepackage{amssymb,amscd,latexsym}   
\usepackage{amsmath}
\usepackage{amsthm}
\newcommand{\llar}{-\kern-5pt-\kern-5pt\longrightarrow}

\newtheorem{Theorem}{Theorem}[section]
\newtheorem{Lemma}[Theorem]{Lemma}
\newtheorem{Corollary}[Theorem]{Corollary}
\newtheorem{Proposition}[Theorem]{Proposition}
\newtheorem{Remark}[Theorem]{Remark}
\newtheorem{Example}[Theorem]{Example}

\newtheorem{Definition}[Theorem]{Definition}

\def\sqr#1#2{{\vcenter{\hrule height.#2pt
        \hbox{\vrule width.#2pt height#1pt \kern#1pt
            \vrule width.#2pt}
        \hrule height.#2pt}}}
\def\phi{\varphi}
\def\demo{\noindent{Proof. }}
\def\square{\mathchoice\sqr64\sqr64\sqr{4}3\sqr{3}3}
\def\qed{\hspace*{\fill} $\square$}

\def\cc{{\bf c}}

\def\xx{{\bf x}}
\def\bb{{\bf b}}
\def\aa{{\bf a}}
\def\yy{{\bf y}}

\def\ff{{\bf f}}

\def\ff{{\bf f}}
\def\gg{{\bf g}}
\def\hh{{\bf h}}
\def\aa{{\bf a}}
\def\qq{{\bf q}}

\def\ii{\'{\i}}

\def\restr{{\kern-1pt\restriction\kern-1pt}}

\def\pp{{\mathbb P}}

\begin{document}
\begin{center}
{\Large{\bf\sc The Newton complementary dual revisited}}

\vspace{0.3in}

{\large\sc Andr\'e D\'oria} \footnote{This author holds a
Post-Doctor Fellowship from the Universidade Federal da Para\ii
ba.} \quad and \quad
 {\large\sc Aron Simis}
\footnote{This author is partially supported by a CNPq grant
(302298/2014-2) and a PVNS Fellowship from CAPES
(5742201241/2016). He thanks the Departamento de Matem\'atica of
the Universidade Federal da Para\ii ba for providing an
appropriate environment for discussions on this work.

D\'oria address: Departamento de Matem\'atica, Universidade Federal
de Sergipe, 49100-000 S\~ao Cristov\~ao, Sergipe, Brazil, email:
avsdoria@gmail.com.

Simis address: Departamento de Matem\'atica, Universidade Federal
de Pernambuco, 50740-640 Recife, PE, Brazil, email:
aron@dmat.ufpe.br.}

\end{center}


\begin{abstract}
	This work deals with the notion of Newton complementary duality as raised originally in the work of the second author and B. Costa. A conceptual revision of the main steps of the notion is accomplished which then leads to a vast simplification and improvement of several statements concerning rational maps and their images. A ring-homomorphism like map is introduced that allows for a close comparison between the respective graphs of a rational map and its Newton dual counterpart.
\end{abstract}

\section{Introduction}
\label{S1}

The notion of a Newton complementary dual to a set of forms having the same degree has been introduced in \cite{NewCremona}.
In this work we improve on some of the aspects and results of the basic theory as exposed in the latter paper, adding a striking unifying simplification to it.
In addition, we apply this simplified approach to settle the behavior of taking Newton dual on notable rational maps and some algebras associated thereof.

In Section 2 we state a brief discussion of the notion and its basic properties, whereas the main core of the section brings up a pack of new relations under the action of the duality. One new fundamental step thereon is the crucial role exerted by the Magnus reciprocal involution
$$(x_0:\cdots :x_n) \mapsto (1/x_0:\cdots :1/x_n)$$
which is the Newton dual to the identity map of $\pp^n$.
A useful observation, not entirely perceived in  \cite{NewCremona} is that taking the Newton dual is nearly as good as evaluating the original forms on the Magnus involution. This in turn yields a uniform procedure to attain other relations.
It allows, e.g., for a short conceptual proof of \cite[Propostion 2.9]{NewCremona} concerning the behavior of rational maps and their images under the action of the Newton duality.
In addition, it has impact on how taking the Newton dual interacts with group-theoretical aspects of Cremona maps.

The main results of the section are Proposition~\ref{PropChapeu1}, Theorem~\ref{TheoChapeuN} and its corollaries.

In Section 3 we develop ab initio the effect of taking the Newton dual on sets of forms defining notable algebras and rational maps.
Here, inspired by the simple behavior of the action of Newton duality on rational maps and their images as obtained in Section 2, we hoped for a parallel behavior at the level of the graph of a rational map, a more conspicuous object, here taken in the disguise of a Rees algebra.
The basic obstruction, enticed by the fact that Newton duality does not induce a ring homomorphism at the level of the polynomial ambient of the Rees algebra, has lead us to introduce a map as a replacement. This map allows to deduce most expected properties of passage between the Rees algebras of a set of forms of the same degree and its Newton dual, in terms of the respective presentation ideals.
Although we do not recover one presentation ideal from the other, we come close in terms of minimal prime ideals thereof.

The other part of the section deals with the interaction of the Newton duality with special elements of the Cremona group.
We chose to work with a subgroup of the Cremona group introduced in \cite{PS2015} based on the notion of de Jonqui\`eres maps. In detail, this group generalizes to higher space the subgroup of plane de Jonqui\`eres maps through a fixed point.
It turns out that taking Newton dual on Cremona maps (i.e., on the corresponding representatives by forms) preserves this subgroup. Moreover, one establishes in terms of the representatives a condition for such a de Jonqui\`eres map to have its Newton dual commute with its inverse map.
Unfortunately, the condition works only in the case where the  de Jonqui\`eres map has a {\em monomial} underlying Cremona map.

The main results of the section are Theorem~\ref{MainRees}, Proposition~\ref{PropDeJonq} and Proposition~\ref{JonqComM}.

\section{General theory}

\subsection{Recap of the Newton complementary dual}

Let $R:= k[\xx] = k[x_0,\ldots,x_n]$ denote the standard graded polynomial ring in $n+1$ indeterminates with coefficients in a field $k$.

Let $f$ be a $d$-form in $R$, with $d\geq 1$.
Denote by $N(f)$ (as a reminder of the Newton polygon) the so-called log-matrix of the (finitely many) terms of $f$ having nonzero coefficients in a fixed ordering, i.e.,
the matrix whose columns are the
exponents vectors of the nonzero terms of $f$ in a fixed ordering.
One may call $N(f)$ the {\sc Newton log matrix} (or simply the {\sc Newton matrix}) of $f$.

Given a finite ordered set $\mathbf{f}:=\{f_0,\ldots,f_m\} (m\geq 1)$ of such forms of the same degree $d\geq 1$,
let $N(\mathbf{f})$ denote the concatenation of the Newton matrices $N(f_0),\ldots,N(f_m)$;
accordingly, we call $N(\mathbf{f})$ the {\sc Newton matrix} of the set $\mathbf{f}$.
Note that $N(\mathbf{f})$ is an integer stochastic matrix.

The row vector $\mathbf{c}_f$ whose entries are the nonzero coefficients of a form $f$ in a fixed ordering is
the {\em coefficient frame} of $f$.
We write symbolically
\begin{equation}\label{symb}
f=\langle\mathbf{c}_f,\xx^{N(f)}\rangle
\end{equation}
as the inner product of the coefficient frame by the set of the corresponding monomials -- it will be called the {\sc Newton representation} of $f$.

The {\sc Newton complementary dual matrix} (or simply the {\sc Newton dual matrix}) of
the Newton matrix $N(\ff)=(a_{i,\ell})$ is the matrix
$$\widehat{N(\ff)}=(\alpha_i - a_{i,\ell}),$$
where $\alpha_i = \max_\ell \{a_{i,\ell}\}$, with $0\leq i\leq n$
and $\ell$ indexes the set of all nonzero terms of all forms in
the set $\mathbf{f}$.

In other words, denoting
${\boldsymbol\alpha}:=(\alpha_0,\ldots,\alpha_n)^t$, one has
$$
\widehat{N(\ff)}=\left[\,{\boldsymbol\alpha}\,| \cdots
|\,{\boldsymbol\alpha}\,\right]_{(n+1) \times (r_0+\ldots + r_m)}
- N(\ff),
$$
where $r_j$ denotes the number of nonzero terms of $f_j$, $j=0,\ldots, m$.
The vector ${\boldsymbol\alpha}$ is called the {\sc directrix vector} of $N(\ff)$ (or of $\ff$ by abuse).
Note that this notion of directrix vector can de defined for any matrix with integer entries and will sometimes employed as such.

For every $j=0,\ldots,m$, let $\widehat{N(\ff)}_j$ denote the submatrix of $\widehat{N(\ff)}$
whose columns come from $f_{j}$.
Finally consider the set of forms defined in terms of their Newton representations
\begin{equation}\label{dual_forms}
\widehat{\ff}:=\{\widehat{f_0}:=\langle\mathbf{c}_0,\xx^{\widehat{N(\ff)}_0}\rangle\,,\ldots, \,\widehat{f_m}:=
\langle\mathbf{c}_m,\xx^{\widehat{N(\ff)}_m}\rangle\}
\end{equation}
where $\mathbf{c}_j=\mathbf{c}_{f_j}$ stands for the coefficient frame of $f_j$.

We call $\widehat{\ff}$ the {\sc Newton complementary dual} set of $\mathbf{f}$.

The reason for calling the procedure a {\em dual} is the following result (\cite[Lemma 2.4]{NewCremona}): if the forms $\ff$ have trivial $\gcd$ and every variable in $\xx$ is in the support of at least one nonzero term of one of the forms, then the Newton dual set $\widehat{\ff}$ also satisfies these conditions and moreover
$$\widehat{\widehat{\ff}}=\ff.$$
These two assumptions were called {\em canonical restrictions} in \cite{NewCremona}.

\smallskip

Let us quote a permeating example

\begin{Example}\label{Magnus}\rm (\cite[Example 2.5]{NewCremona}) If $\xx=\{x_0,\ldots, x_n\}$ then its directrix vector is $(1,\ldots,1)^t$ and
    \begin{eqnarray}\label{Magnus_involution}
     \widehat{\xx}=\{x_1\cdots x_n,\, \ldots, \, x_0\cdots \underline{x_i}\cdots x_n,
    \,\ldots, \,x_0\cdots x_{n-1}\},
    \end{eqnarray}
    where $\underline{x_i}$ means $x_i$ is omitted.
    In terms of the rational maps defined by the respective forms, this says that the Newton complementary dual of the identity
    map of $\pp^n$ is the Magnus reciprocal involution.
\end{Example}

This simple example will be quite crucial in expressing basic relations as in the next subsection.

\subsection{New relations}

In this section we state a few relations satisfied by taking the Newton dual, some of which greatly simplify some parts of \cite{NewCremona}.
We start by observing the impact of evaluating forms on the Newton dual $\widehat{\xx}$ of the set $\xx$ of the variables.

\begin{Lemma}\label{LemNewtComp}
Let $g \in R := k[x_0, \ldots , x_n]$ be arbitrary form and
$\widehat{\xx}$ denote the Newton complementary dual set of $\xx
:= \{x_0,\ldots,x_n\}$. Then $N(g(\widehat{\xx}))=N(\widehat{\xx}).N(g)$,
where the dot indicates matrix multiplication.
In particular, the Newton representation of $g(\widehat{\xx})$ is
$$g(\widehat{\xx}) = \langle \cc_g , \xx^{N(\widehat{\xx}).N(g)} \rangle.$$
\end{Lemma}
\demo Let $g_1,g_2\in R$ be forms with empty common monomial support (i.e., no nonzero term appears in both forms) satisfying the desired equality in the statement. Then, up to suitable ordering of the constituent monomials, one has
\begin{eqnarray*}
N((g_1+g_2)(\widehat{\xx}))&=&N(g_1(\widehat{\xx})+g_2(\widehat{\xx}))=N(g_1(\widehat{\xx}))\sqcup N(g_2(\widehat{\xx}))\\
&=& N(\widehat{\xx}).N(g_1) \sqcup N(\widehat{\xx}).N(g_2)=N(\widehat{\xx}).(N(g_1) \sqcup N(g_2))\\
& = & N(\widehat{\xx}).N(g_1+g_2),
\end{eqnarray*}
where $\sqcup$ denotes matrix concatenation.

The supplementary result about the Newton representation is clear since, in the assumed situation of the forms $g_1,g_2$, one has the relation $\mathbf{c}_{g_1+g_2}=\mathbf{c}_{g_1}\sqcup \mathbf{c}_{g_2}$ for the respective coefficient frames.

Therefore we have reduced the problem to the case when $g(\xx) = c
\: x_0^{a_0}\ldots \: x_n^{a_n}$ is a term. In this situation the result is obvious because $$N(g) = \left(%
\begin{array}{c}
   a_0 \\
   \vdots\\
   a_n\\
\end{array}%
\right) \; {\rm and}\; N(\widehat{\xx})= \left(\begin{array}{ccccc}
                                  0 & 1 & 1 & \cdots & 1 \\
                                  1 & 0 & 1 & \cdots & 1 \\
                                  1 & 1 & 0 & \cdots & 1 \\
                                  \vdots & \vdots & \vdots & \ddots & \vdots \\
                                  1 & 1 & 1 & \cdots & 0 \\
                                \end{array}\right).$$
\qed

\begin{Proposition}\label{PropChapeu1}
Let $g \in k[x_0, \ldots , x_n]$ denote an arbitrary form of degree $d\geq 1$
and let $\widehat{g}$ denote its Newton complementary dual. Then
$$
g(\widehat{\xx}) = x_0^{d - \beta_0} \cdots \: x_n^{d - \beta_n}
\: \widehat{g},
$$ where  $\widehat{\xx}$ denote the Newton
complementary dual set of $\xx := \{x_0,\ldots,x_n\}$ and
${\boldsymbol\beta} := (\beta_0, \ldots, \beta_n)^t$ is the directrix
vector of $N(g)$.
\end{Proposition}
\demo By Lemma \ref{LemNewtComp},
$$
N(g(\widehat{\xx}) )= N(\widehat{\xx}).N(g).
$$
Therefore, by taking respective Newton representations it suffices to prove the matrix equality
$$
N(\widehat{\xx}).N(g) =
[{\boldsymbol\alpha}|\ldots|{\boldsymbol\alpha}]_{(n+1) \times r}
+ \widehat{N(g)},
$$
where $r$ denotes the number of nonzero terms of $g$ and
${\boldsymbol\alpha} = (d - \beta_0, \ldots, d - \beta_n)^t$.
But this is clear as
$$\begin{array}{ccl}
    N(\widehat{\xx}).N(g) & = & (\mathbb{I} - N(\xx)).N(g) = \mathbb{I} .N(g) - N(\xx).N(g) = \mathbb{I} .N(g) - N(g)\\
                          & = &  \mathbb{I} .N(g) - [\beta | \ldots | \beta] + [\beta | \ldots | \beta] - N(g)\\
                          & = &
                          [{\boldsymbol\alpha}|\ldots|{\boldsymbol\alpha}] +
                          \widehat{N(g)},
\end{array}$$
where $\mathbb{I}$ is the  $(n+1)\times (n+1)$ matrix whose entries are all $1$.
 \qed

\medskip

Thus, the proposition tells us that to get the Newton dual $\widehat{g}$ of a form $g$ one has to evaluate $g$ on the terms of the involution in (\ref{Magnus_involution}) and multiply out by a certain Laurent monomial.
The moral of this statement is that it may be easier to handle the properties of  $\widehat{g}$ by evaluation than by applying the original definition.
Thus, for example, one can derive the following result:

\begin{Corollary}\label{ComNewMult}
Let $p,q\in R$. Then $\widehat{pq} = \widehat{p}\, \widehat{q}.$
\end{Corollary}
\demo Set $g := p.q$. Clearly, $g(\widehat{\xx}) =
p(\widehat{\xx}).q(\widehat{\xx}),$ where $\widehat{\xx}$ denote
the Newton complementary dual set of $\xx = \{x_0,\ldots,x_n\}$.
By Proposition~\ref{PropChapeu1}, one has
$$g(\widehat{\xx}) = x_0^{r+s -
\gamma_0} \cdots \: x_n^{r+s - \gamma_n} \: \widehat{g}, \:
p(\widehat{\xx}) = x_0^{r - \alpha_0} \cdots \: x_n^{r - \alpha_n}
\: \widehat{p},\, q(\widehat{\xx}) = x_0^{s - \beta_0}
\cdots \: x_n^{s - \beta_n} \: \widehat{q},$$
where $r=\deg(p), s=\deg(q)$ and
${\boldsymbol\gamma} := (\gamma_0, \ldots, \gamma_n)^t$,
${\boldsymbol\alpha} := (\alpha_0, \ldots, \alpha_n)^t$ and
${\boldsymbol\beta} := (\beta_0, \ldots, \beta_n)^t$ are the
directrix vectors of $N(g)$, $N(p)$ and $N(q)$, respectively.
Thus
$$
x_0^{r+s - \gamma_0} \cdots \: x_n^{r+s - \gamma_n} \: \widehat{g}
= x_0^{r+s - \alpha_0 - \beta_0} \cdots \: x_n^{r+s - \alpha_n -
\beta_n} \: \widehat{p}.\widehat{q}.
$$
On the other hand, one has $\deg_{x_i}(g) = \deg_{x_i}(p) + \deg_{x_i}(q)$, for every $i = 0, \ldots,n$.
Therefore ${\boldsymbol\gamma} = {\boldsymbol\alpha} +
{\boldsymbol\beta}$ and hence we are through.\qed

\medskip

 Proposition~\ref{PropChapeu1} extends to a set of forms:

\begin{Theorem}\label{TheoChapeuN}
Let $\gg := \{g_0,\ldots ,g_m\} \subset R := k[x_0, . . . , x_n]$
denote a set of arbitrary forms of same degree $d\geq 1$. Then
$$\gg(\widehat{\xx}) = x_0^{d - \beta_0} \cdots \: x_n^{d -
\beta_n} \: \widehat{\gg},$$ where $\widehat{\xx}$ denote the
Newton complementary dual set of $\xx := \{x_0,\ldots,x_n\}$ and
${\boldsymbol\beta} := (\beta_0, \ldots, \beta_n)^t$ is the directrix
vector of $N(\gg)$.
\end{Theorem}
\demo Applying Proposition~\ref{PropChapeu1} with $g:=g_i$ gives
$$g_i(\widehat{\xx}) = x_0^{d - \beta_{0,i}} \cdots \:
x_n^{d - \beta_{n,i}} \: \widehat{g_i},$$ where
$\underline{\beta_i} := (\beta_{0,i}, \ldots, \beta_{n,i})^t$
 is the directrix vector of $g_i$.

By definition, the directrix vector ${\boldsymbol\beta} =
(\beta_0, \ldots, \beta_n)^t$  of $\gg = \{g_0,\ldots,g_m\}$
satisfies the inequalities $\beta_i - \beta_{i,j} \geq 0$, for all
$i = 0,\ldots,n$ and $j = 0,\ldots,m$. Therefore
$$\begin{array}{ccl}
  x_0^{d - \beta_{0,i}} \cdots \: x_n^{d - \beta_{n,i}} \: \widehat{g_i} & = &  x_0^{d - \beta_{0,i}} \cdots \: x_n^{d - \beta_{n,i}} \: \displaystyle\frac{x_0^{\beta_0 - \beta_{0,i}} \ldots \: x_n^{\beta_n - \beta_{n,i}}}{x_0^{\beta_0 - \beta_{0,i}} \ldots \: x_n^{\beta_n - \beta_{n,i}}} \: \widehat{g_i}\\
  & = &  x_0^{d - \beta_0} \cdots \: x_n^{d - \beta_n} \: x_0^{\beta_0 - \beta_{0,i}} \cdots \: x_n^{\beta_n - \beta_{n,i}} \: \widehat{g_i}.\\
\end{array}$$
Note that
$$
N(\widehat{\gg}) = \left( [ \gamma_0 | \ldots |\gamma_0]_{(n+1)
\times r_0} + N(\widehat{g_0}) \:\: \ldots \:\: [ \gamma_m |
\ldots |\gamma_m]_{(n+1) \times r_m} + N(\widehat{g_m}) \right),
$$ where $\gamma_j := (\beta_0 -
\beta_{0,j},\ldots,\beta_n - \beta_{n,j})^t$ and $r_j$ denotes the
number of nonzero terms of $g_j,$ $j = 0,\ldots, m$. Therefore,
$\gg(\widehat{\xx}) = x_0^{d - \beta_0} \cdots \: x_n^{d -
\beta_n} \: \widehat{\gg}$, as was to be shown. \qed

\medskip

Next is a generalization of the previous proposition to the case where one evaluates on the forms of the Newton dual of a set of $n+1=\dim R$ forms of the same degree (instead of the variables $\xx$).

\begin{Proposition}\label{PropCompChapeuN}
Let $\gg := \{g_0,\ldots ,g_m\} \subset R$ be arbitrary forms of degree $d\geq 1$ and let $\hh := \{h_0,\ldots ,h_n\} \subset R$
denote a set of $n+1=\dim R$ forms of a fixed degree.
Then
$$\gg(\widehat{\hh}) = x_0^{d \alpha_0 - \beta_0} \cdots \: x_n^{d
\alpha_n - \beta_n} \: \widehat{\gg(\hh)},$$ where
${\boldsymbol\alpha} := (\alpha_0, \ldots, \alpha_n)^t$ and
${\boldsymbol\beta} := (\beta_0, \ldots, \beta_n)^t$ are the
directrix vectors of $N(\hh)$ and $N(\gg(\hh))$, respectively.
\end{Proposition}
\demo We apply  Theorem~\ref{TheoChapeuN} twice, first taking $\gg$ to be $\hh$; it obtains
$$
\hh(\widehat{\xx}) = x_0^{s - \alpha_0} \cdots \: x_n^{s -
\alpha_n} \: \widehat{\hh},
$$
where $s$ is the common degree of the forms $\hh$. Evaluating $\gg$ on these $n$ forms yields
$$
\begin{array}{ccl}
  \gg(\widehat{\hh}) & = & (x_0^{ds - d\alpha_0} \cdots \:
x_n^{ds - d\alpha_n})^{-1} \: \gg(\hh(\widehat{\xx})) \\
   & = & (x_0^{ds - d\alpha_0} \cdots \:
   x_n^{ds - d\alpha_n})^{-1}  \:
(\gg(\hh))(\widehat{\xx}).\\
\end{array}
$$ Apply Theorem~\ref{TheoChapeuN} again, this time around taking $\gg$ to be  $\gg(\hh)$. It now yields
$$
\gg(\hh)(\widehat{\xx}) = x_0^{ds - \beta_0} \cdots \: x_n^{ds -
\beta_n} \: \widehat{\gg(\hh)}.
$$ Therefore
$$
\gg(\widehat{\hh}) = (x_0^{ds - d\alpha_0} \cdots \:
x_n^{ds - d\alpha_n})^{-1}\: x_0^{ds - \beta_0} \cdots \:
x_n^{ds - \beta_n} \: \widehat{\gg(\hh)} = x_0^{d \alpha_0 -
\beta_0} \cdots \: x_n^{d \alpha_n - \beta_n} \:
\widehat{\gg(\hh)}.
$$\qed

\medskip

\begin{Corollary}\label{isomofphism}
    Let $\gg := \{g_0,\ldots,g_m\} \subset R$ be a set of forms of the
    same degree satisfying the canonical restrictions. Then $k[\gg]$
    and $ k[\widehat{\gg}]$ are isomorphic as $k$-algebras.
\end{Corollary}
\demo It suffices to show that any homogeneous polynomial relation
of $\gg$ is one of $\widehat{\gg}$, and vice-versa. Let
$F(\yy) \in k[\yy] := k[y_0,\ldots,y_m]$ be a homogeneous
polynomial of degree $d$ such that $F(\gg)=0$.
Since $\gg$ satisfies the canonical restrictions, applying Theorem~\ref{TheoChapeuN} to the Newton dual $\widehat{\gg}$ yields

$$\widehat{\gg}(\widehat{\xx})=M\, \widehat{\widehat{\gg}}=M\,\gg,
$$
 for suitable monomial $M \in R$.
 Therefore,
 $$F(\widehat{\gg}(\widehat{\xx}))=F(M\,\gg)=M^d\, F(\gg)=0.$$

But then $F(\widehat{\gg}(\widehat{\xx}))=F(\widehat{\gg})(\widehat{\xx})$  tells us that the form $F(\widehat{\gg})$ further evaluated at the forms of $\widehat{\xx}$ vanishes.
Since the latter forms are algebraically independent over $k$ it follows that  $F(\widehat{\gg})=0$.

The other direction is similar.
\qed

\medskip

For the next corollary, we note that if $\mathfrak{G}: \pp^n\dasharrow \pp^m$ is a rational map defined by a set of forms $\gg = \{g_0,\ldots, g_m\} \subset R$ of the same degree then the rational map obtained by composing the Magnus involution with $\mathfrak{G}$ is defined by the evaluated forms $\gg(\widehat{\xx}) = \{g_0(\widehat{\xx}),\ldots, g_m(\widehat{\xx})\}$.
This observation will be used next without further ado.

\begin{Corollary}\label{corobira}
    Let $\gg = \{g_0,\ldots, g_m\} \subset R := k[x_0,\ldots, x_n] (m\geq 1, n
    \geq 1)$ be arbitrary forms of the same degree and let
    $\widehat{\gg}$ denote its Newton complementary dual set. Then:
    \begin{itemize}
        \item[(a)] $\gg(\widehat{\xx})$ and $\widehat{\gg}$ define the same rational
        map.
        \item[(b)] If $\gg$ defines a birational map onto its image then so does $\widehat{\gg}$ and the two maps have the same image.
    \end{itemize}
\end{Corollary}
\demo (a) This follows immediately from Theorem~\ref{TheoChapeuN} provided, as is common practice, we do not distinguish between two rational maps whose defining tuples match up to multiplication by a nonzero element of the fraction field $k(x_0,\ldots, x_n)$ of $R$.

(b) This follows from (a)  and Corollary~\ref{isomofphism}.
\qed

\begin{Corollary}\label{compMagnus}
    Let $\mathfrak{G}: \pp^n\dasharrow \pp^n$ denote a Cremona map, let $\widehat{\mathfrak{G}}$ denote the corresponding Newton dual Cremona map obtained via {\rm Corollary~\ref{corobira}(b)},
     and let $\mathfrak{M}$ stand for the Magnus
    reciprocal involution. Then, in terms of the group law in the Cremona group of $\pp^n$, one has:
    \begin{itemize}
        \item[(i)] $\widehat{\mathfrak{G}} = \mathfrak{G} \circ \mathfrak{M}$ and $\left(\widehat{\mathfrak{G}}\right)^{-1} =
        \mathfrak{M} \circ \mathfrak{G}^{-1}.$
        \item[(ii)] $\widehat{\mathfrak{G}^{-1}} = \left(\widehat{\mathfrak{G}}\right)^{-1}$
        if and only if $\mathfrak{M} \circ \mathfrak{G} = \mathfrak{G} \circ
        \mathfrak{M}.$

        In particular, if $\mathfrak{G}$ is a monomial Cremona map, then $\widehat{\mathfrak{G}^{-1}} = \left(\widehat{\mathfrak{G}}\right)^{-1}$.
        \item[(iii)] Let $C,D \in k[\xx]$ denote the source inversion factors of $\mathfrak{G}$ and $\widehat{\mathfrak{G}}$, respectively, one has: $$\Gamma_1 \left(\widehat{C}\right)^n = q_1 D \:\: \mbox{ and } \:\: \Gamma_2 \left(\widehat{D}\right)^n = q_2 C,$$ where $\Gamma_1,\Gamma_2,q_1,q_2 \in k[\xx]$, with $\Gamma_1$ and $\Gamma_2$ are monomials.
    \end{itemize}

\end{Corollary}
\demo (i) This is a rephrasing of Corollary \ref{corobira} (a).

(ii) By item (i),
$\left(\widehat{\mathfrak{G}}\right)^{-1} = \mathfrak{M} \circ
\mathfrak{G}^{-1}$. By the same token, $\widehat{\mathfrak{G}^{-1}}
= \mathfrak{G}^{-1} \circ \mathfrak{M}$. Therefore,
$$\widehat{\mathfrak{G}^{-1}} =
\left(\widehat{\mathfrak{G}}\right)^{-1} \: \Leftrightarrow \:
\mathfrak{G}^{-1} \circ \mathfrak{M} = \mathfrak{M} \circ
\mathfrak{G}^{-1} \: \Leftrightarrow \: \mathfrak{M} \circ
\mathfrak{G} = \mathfrak{G} \circ \mathfrak{M}.$$

The supplementary assertion follows from the fact that $\mathfrak{G}$
commutes with $\mathfrak{M}$, when $\mathfrak{G}$ is a monomial
Cremona map.

(iii) Let $\gg = \{g_0,\ldots,g_n\}$ and $\gg' = \{g'_0,\ldots,g'_n\}$ denote the representative of $\mathfrak{G}$ and $\mathfrak{G}^{-1}$, respectively, satisfying ${\rm gcd}\{g_0,\ldots, g_n\} = {\rm gcd}\{g'_0,\ldots, g'_n\} = 1$. Thus
$$g'_i(g_0,\ldots,g_n ) = x_i \: C,$$ for every $i = 0,\ldots,n$. By Theorem~\ref{TheoChapeuN}, $\widehat{\gg} = \frac{1}{M} \gg(\widehat{\xx})$, where $M$ is a suitable monomial. Let $\hh = \{h_0,\ldots,h_n\}$ denote the representative of $\left(\widehat{\mathfrak{G}}\right)^{-1}$ satisfying ${\rm gcd}\{h_0,\ldots, h_n\} = 1$. So $$h_i(\widehat{\gg}) = x_i \: D,$$ for every $i = 0,\ldots,n$. By (i), $p \: \hh = \widehat{\xx}(\gg')$, for suitable form $p$.

Therefore, for a $s$ non-negative integer,
$$
\begin{array}{cl}
  x_0 \: D & = h_0(\widehat{\gg}) = \left(\frac{1}{p} \: g'_1 \cdots g'_n \right) \circ \left( \frac{1}{M} \gg(\widehat{\xx}) \right) =  \frac{1}{M^s} \left[\left(\frac{1}{p} \: g'_1 \cdots g'_n \right) \circ \left( \gg(\widehat{\xx})\right)\right]\\
  & = \frac{1}{M^s \: p(\gg(\widehat{\xx}))} \left[\left(g'_1 \cdots g'_n \right) \circ \left( \gg(\widehat{\xx})\right)\right] = \frac{1}{M^s \: p(\gg(\widehat{\xx}))} \left[\left(g'_1(\gg) \cdots g'_n(\gg) \right) \circ \widehat{\xx}\right]\\
  & = \frac{1}{M^s \: p(\gg(\widehat{\xx}))} \left[\left(x_1 \cdots x_n C^n \right) \circ \widehat{\xx}\right] = \frac{x^n_0 x^{n-1}_1 \ldots \: x_n^{n-1}}{M^s \: p(\gg(\widehat{\xx}))}
  [C(\widehat{\xx})]^n.\\
\end{array}
$$ By Proposition~\ref{PropChapeu1}, with $g := C$, $(\beta_0, \ldots, \beta_n)^t$ is the directrix vector of $C$ and $d := {\rm deg}(C)$, $$\Gamma_2 \left(\widehat{C}\right)^n = q_2 D,$$ where $\Gamma_2 := x^{n-1 + nd - n\beta_0}_0 \ldots \: x_n^{n-1 + nd - n\beta_n}$ and $q_2 := M^s \: p(\gg(\widehat{\xx}))$. The other equation in (iii) following the duality of the Newton complementary dual.\qed

\bigskip

We have seen that taking Newton dual on individual forms is a multiplicative operation.
Additivity does not hold, but there is a close relationship. Moreover, one can actually deal with sets in the following sense: given sets of forms $\mathbf{p} = \{p_0,\ldots,p_m\} \subset R$ and $\qq =\{q_0,\ldots,q_m\}\subset R$, both in the same degree, we write $\mathbf{p}+\qq$ for the set of forms obtained by adding the tuples $(p_0,\ldots,p_m)$ and $(q_0,\ldots,q_m)$ coordinate wise.

\begin{Corollary}
Let $\mathbf{p} = \{p_0,\ldots,p_m\} \subset R$ and $\qq =
\{q_0,\ldots,q_m\}$ be sets of forms of $R$, both in the same degree
$d$.
Then
$$( x_0^{d - \gamma_0} \cdots \: x_n^{d - \gamma_n}) \: \widehat{\mathbf{p} + \qq} =
    (x_0^{d - \alpha_0} \cdots \: x_n^{d - \alpha_n}) \: \widehat{\mathbf{p}} + (x_0^{d - \beta_0} \cdots \: x_n^{d - \beta_n}) \: \widehat{\qq},$$
where
    ${\boldsymbol\gamma} := (\gamma_0, \ldots, \gamma_n)^t$,
    ${\boldsymbol\alpha} := (\alpha_0, \ldots, \alpha_n)^t$ and
    ${\boldsymbol\beta} := (\beta_0, \ldots, \beta_n)^t$ are the
    directrix vectors of $\mathbf{p} + \qq$, $\mathbf{p}$ and $\qq$, respectively.
\end{Corollary}
\demo This is an immediate application of Theorem~\ref{TheoChapeuN} to all sets (tuples) of the relation $(\mathbf{p}+\qq)(\widehat{\xx})=\mathbf{p}(\widehat{\xx}) + \qq(\widehat{\xx})$.
\qed

\section{Interaction at large}

\subsection{Interaction with the graph of a rational map}

Let $\gg = \{g_0,\ldots,g_m\} \subset R:= k[\xx] =
k[x_0,\ldots,x_n]$ be arbitrary forms of the same degree. In this
section we will see a relationship between the defining ideals of
the Rees algebras of $\gg$ and $\widehat{\gg}$. Our reference for
Rees algebra is \cite{bir2003}, which contains the material
in the form we use here.

Let $R[\yy] := R[y_0,\ldots, y_m]$ denote a polynomial ring over R
with the standard bigrading where $\deg(x_i) = (1,0)$ and
$\deg(y_j) = (0,1)$.

Let $p=p(\xx,\yy) \in R[\yy]=k[\xx,\yy]$ denote a bihomogeneous
polynomial of bidegree $(d_x,d_y)$. Write

$$p(\xx,\yy) = \sum_{\stackrel{ |\aa|=d_x}{\stackrel{|\bb|=d_y}{}}} \: c_{\aa, \bb}\: \xx^{\aa} \yy^{\bb},$$
where $\aa=(a_0,\ldots,a_n), \bb=(b_0,\ldots,b_m)$.

As in the previous part, we let $\cc(p)$ denote the row of nonzero
coefficients $ c_{\aa, \bb}$ of $p$ in a fixed ordering (the coefficient frame of
$p$) and let $N_{\xx}(p)$ (respectively, $N_{\yy}(p)$) denote the
(Newton like) matrix whose columns are the exponent vectors of the
corresponding $\xx$-terms (respectively, $\yy$-terms) in the same
ordering.

As before, we may write symbolically
\begin{equation}\label{symbol_form}
p(\xx,\yy)= \langle \cc(p),\xx^{N_{\xx}(p)} \yy^{N_{\yy}(p)}
\rangle
\end{equation}
and call it the {\em Newton representation} of $p$ in analogy to the homogeneous case.

The following result is an analogue of the
Proposition~\ref{PropChapeu1} for bihomogeneous polynomials.

\begin{Lemma}\label{LemmDualbiho}
Let $p(\xx,\yy) \in R[\yy]$ be a bihomogeneous form such that $d_x\geq 1$. Then,
$$p(\widehat{\xx},\yy) = x_0^{d_x - \beta_0} \cdots \: x_n^{d_x - \beta_n} \: \langle \cc(p) , \xx^{\widehat{N_\xx(p)}}\yy^{N_\yy(p)}
\rangle,$$ where ${\boldsymbol\beta} := (\beta_0, \ldots,
\beta_n)^t$ is directrix vector of $N_\xx(p)$ and $\widehat{\xx}$
denotes the Newton complementary dual set of $\xx =
\{x_0,\ldots,x_n\}$.
\end{Lemma}
\demo
For the sake of clarity, denote $P(\xx):=p(\xx,\yy)$ as an element of $K[\xx]$, where $K=k(\yy)$ stands for the field of fractions of $k[\yy]$.
Applying Proposition~\ref{PropChapeu1} to $P(\xx)$ gives
\begin{equation}\label{BigP}
P(\widehat{\xx})= (x_0^{d_x - \beta_0} \cdots \: x_n^{d_x - \beta_n} ) \widehat{P(\xx)},
\end{equation}
since $N(P(\xx))=N_x(p(\xx,\yy))=N_x(p)$ as one easily sees from the definitions, and hence the directrix vector does not change.
At the other end, one has the Newton representation of the Newton dual of $P(\xx)$:
\begin{equation}\label{Newton_BigP}
\widehat{P(\xx)}=\langle \mathbf{c}(\widehat{P(\xx)}), \xx^{N_x(p)}\rangle,
\end{equation}
where $\mathbf{c}(\widehat{P(\xx)})=\mathbf{c}(P(\xx))$ denotes the coefficient frame of $P(\xx)$ over $K$.
The elements of this frame (in the established ordering) are the terms of the form $c_{\aa,\bb}\yy^{\bb}$.
Therefore, taking in account the definition of the Newton matrix $N_y(p)$, the required expression follows from  (\ref{BigP}) and (\ref{Newton_BigP}).
\qed

\medskip

The basic result of this part is the following:

\begin{Proposition}\label{evaluating}
    Let $\gg = \{g_0,\ldots,g_m\} \subset R$ stand for a set of forms of the same degree and let $p(\xx,\yy)\in R[\yy]$ denote a bihomogeneous form.
    If $p(\xx,\gg)=0$ then $p(\widehat{\xx},\widehat{\gg})=0$.
\end{Proposition}
\demo
By Theorem~\ref{TheoChapeuN}, there exists a suitable monomial $M$ in $k[\xx]$ such that
$$
p(\widehat{\xx},\widehat{\gg}) = p\left(\widehat{\xx},\frac{\gg(\widehat{\xx})}{M}\right) = \frac{1}{M^{d_y}}\, p(\widehat{\xx},\gg(\widehat{\xx}))=\frac{1}{M^{d_y}}\, p(\xx,\gg)(\widehat{\xx}).
$$
Therefore, $p(\widehat{\xx},\widehat{\gg}) =0$, as was to be shown.
\qed

\begin{Remark}\rm Since we have only defined the Newton dual for forms of positive degree, both Proposition~\ref{PropChapeu1} and Lemma~\ref{LemmDualbiho} require this assumption.
However, one can establish the  convention that for a form $p$ of degree zero,
$N_\xx(p)$ is the zero matrix.
In this way, the above lemma also works in the case $d_x=0$, in which case one recovers the Newton representation of a form in $k[\yy]$.
\end{Remark}

With this convention and the notation of Lemma~\ref{LemmDualbiho} we introduce the following map of the standard graded polynomial ring $k[\xx,\yy]$:
\begin{Definition}\rm
    Let $p(\xx,\yy) \in R[\yy]$ be a bihomogeneous form. Set
    $$\psi(p) := \langle \cc(p) ,
    \xx^{\widehat{N_\xx(p)}}\yy^{N_\yy(p)} \rangle = ( x_0^{d_x - \beta_0} \cdots \: x_n^{d_x - \beta_n})^{-1}\, p(\widehat{\xx},\yy).$$
\end{Definition}
The definition of $\psi$ is naturally extended to the whole of $k[\xx,\yy]$ using the direct sum decomposition $k[\xx,\yy]=\bigoplus_{(a,b)\in \mathbb{N}} k[\xx,\yy]_{a,b}$, where $k[\xx,\yy]_{a,b}$ denotes the $k$-vector space spanned by the monomials of bidegree $(a,b)$.

Note that, with the above convention, one has $\psi(q)=q$ for any homogeneous form $q\in k[\yy]$. In other words, the restriction of $\psi$ to the $k$-subalgebra $k[\yy]\subset k[\xx,\yy]$ is the identity map.

Unfortunately, $\psi$ is not a ring homomorphism, but it is quite close in the following manner:

\begin{Lemma}\label{psi+.}
Let $p,q \in k[\xx,\yy]$ bihomogeneous polynomials. One has:
\begin{itemize}
\item[(i)]  $M\psi(p + q) = M_1\psi(p)+M_2\psi(q)$, where $M$, $M_1$ and $M_2$
are suitable monomials in $k[\xx]$.
\item[(ii)] $M\psi(pq) = M' \psi(p) \psi(q)$, where $M$ and $M'$
are suitable monomials in $k[\xx]$.
\end{itemize}
\end{Lemma}
\demo (i) If $p$ and $q$ have different bidegrees then the result is trivial by the definition of $\psi$ extended to the whole of $k[\xx,\yy]$.
Thus, assume they have the same bidegree.
Set $g := p + q$.
Clearly,  $$g(\widehat{\xx},\yy) =
p(\widehat{\xx},\yy) + q(\widehat{\xx},\yy).$$
Then the result follows immediately from Lemma~\ref{LemmDualbiho}.

(ii) This proof is analogous to the one in item (i). \qed

\bigskip

Fix presentations
\begin{equation}\label{Rees_presentation}
\mathcal{R}_R(\gg)\simeq R[\yy]/\mathcal{J}_{\gg} \: \text{\rm and} \:
\mathcal{R}_R(\widehat{\gg})\simeq R[\yy]/\mathcal{J}_{\widehat{\gg}}
\end{equation}
of the respective Rees algebras of $\gg$ and $\widehat{\gg}$.
The ideals $\mathcal{J}_\gg$ and $\mathcal{J}_{\widehat{\gg}}$ are often called {\em presentation ideals}.

The main result of this part is the following:

\begin{Theorem}\label{MainRees}
Let $\gg = \{g_0,\ldots,g_m\} \subset R$ be a set of forms of the
same degree satisfying the canonical restrictions and let
$\mathcal{J}_\gg$ and $\mathcal{J}_{\widehat{\gg}}$ be as in {\rm (}\ref{Rees_presentation}{\rm )}.
Then:
\begin{enumerate}
\item[{\rm (a)}] $\psi$ maps $\mathcal{J}_\gg$ to $\mathcal{J}_{\widehat{\gg}}$.
\item[{\rm (b)}] Given a bihomogeneous $q(\xx,\yy)\in \mathcal{J}_{\widehat{\gg}}$ such that no variable of $k[\xx]$ is a factor of $q(\xx,\yy)$ then $q(\xx,\yy) = \psi(p(\xx,\yy))$ for some bihomogeneous  $p(\xx,\yy) \in \mathcal{J}_\gg$;
 In particular, a finite subset of $\mathcal{J}_\gg$ maps onto a set of minimal generators of  $\mathcal{J}_{\widehat{\gg}}$.
\item[{\rm (c)}] If $\widetilde{\mathcal{J}_\gg}\subset \mathcal{J}_{\widehat{\gg}}$ is the subideal generated by the image by $\psi$ of a set of generators of $\mathcal{J}_\gg$ then $\mathcal{J}_{\widehat{\gg}}$ is a minimal prime of $\widetilde{\mathcal{J}_\gg}$.
\end{enumerate}
\end{Theorem}
\demo (a) This follows straightforwardly from Proposition~\ref{evaluating} and the definition of $\psi$.

(b) Since $\widehat{\widehat{\gg}} = \gg$, it suffices to prove that
$\psi(\psi(q(\xx,\yy))) = q(\xx,\yy)$.  By
hypothesis, no  $x_i$ is factor of $q(\xx,\yy)$, so
every row of the Newton matrix $N_\xx(q)$ has at least one column such that the
$i$th column entry vanishes. Therefore,
$\widehat{\widehat{N_\xx(q)}} = N_\xx(q)$, that is,
$\psi(\psi(q(\xx,\yy))) = q(\xx,\yy)$.

(c) We will show a slightly stronger result, namely, there is a monomial $M\in k[\xx]$ such that $\mathcal{J}_{\widehat{\gg}} \subset
\widetilde{\mathcal{J}_\gg} : M$.
Since $M\notin \mathcal{J}_{\widehat{\gg}}$ because the presentation ideal is prime, then this will force the equality $\mathcal{J}_{\widehat{\gg}} =
\widetilde{\mathcal{J}_\gg} : M$.

Let $q(\xx,\yy)$ denote a bihomogeneous minimal generator of $\mathcal{J}_{\widehat{\gg}}$.
In particular, by part (b) there exists $p(\xx,\yy) \in \mathcal{J}_\gg$ such that $\psi(p(\xx,\yy)) =
q(\xx,\yy)$.
Now express $p(\xx,\yy)$ in terms of the given (finite) set of generators of $\mathcal{J}_\gg$. By Lemma~\ref{psi+.}, one can write
$$ M\psi(p(\xx,\yy)) = \sum_{i = 1}^t h_i \psi(p_i(\xx,\yy))$$
for suitable monomial $M\in k[\xx]$ and certain biforms $\{h_0,\ldots,h_t\}\in k[\xx,\yy]$.
Therefore $M \,q(\xx,\yy) \in \widetilde{\mathcal{J}_\gg}$.
Now, letting $q(\xx,\yy) $ run through a finite set of minimal (hence, irreducible) generators of $\mathcal{J}_{\widehat{\gg}}$, and letting $M$ denote the product of the corresponding finitely many monomials of $k[\xx]$, one obtains
$\mathcal{J}_{\widehat{\gg}} \subset \widetilde{\mathcal{J}_\gg}
: M$, as promised.
\qed

\smallskip

\begin{Example}\rm
Let $\gg = \{x_1^2x_3^3, x_0x_1^3x_2, x_1x_2^3x_3, x_2^3x_3^2\}\subset R=k[x_0,x_1,x_2,x_3]$.
\end{Example}

A calculation with \cite{singular} gives that $\mathcal{J}_\gg$  minimally generated by the biforms
\begin{eqnarray*}\nonumber
x_3y_2-x_1y_3,\;x_2^3y_0-x_1^2x_3y_3,\;
x_0x_1x_2y_0-x_3^3y_1,\;
x_2^2x_3y_1-x_0x_1^2y_2, \\  \nonumber
 x_0x_1y_2^2-x_2^2y_1y_3, \;x_0x_2y_0y_2-x_3^2y_1y_3,\; x_0^2y_0y_2^4-x_2x_3y_1^2y_3^3.
\end{eqnarray*}
Similarly, a computation with $\widehat{\gg} = \{x_0x_1x_2^3, x_2^2x_3^3, x_0x_1^2x_3^2, x_0x_1^3x_3\}$ gives the following set of minimal generators of
$\mathcal{J}_{\widehat{\gg}}$:
\begin{eqnarray*}\nonumber
x_1y_2-x_3y_3,\;x_0x_1^2y_1-x_2^2x_3y_2,\;  x_1^2x_3y_0-x_2^3y_3,\; x_3^3y_0-x_0x_1x_2y_1, \; x_1x_3^2y_0-x_2^3y_2,\\
x_2^2y_2^2-x_0x_1y_1y_3,\; x_3^2y_0y_2-x_0x_2y_1y_3,\;
 x_2x_3y_0y_2^4-x_0^2y_1^2y_3^3.
\end{eqnarray*}
Note the generator $x_1x_3^2y_0-x_2^3y_2$ of $\mathcal{J}_{\widehat{\gg}}$ is the image of the form $x_2^3y_0-x_1x_3^2y_2$ by $\psi$.
The latter is not a member of the above minimal set of generators of  $\mathcal{J}_\gg$.
However, note that $x_1x_3(x_1x_3^2y_0-x_2^3y_2) \in
\widetilde{\mathcal{J}_\gg}.$
A direct verification gives that the remaining members of the above set of generators of of $\mathcal{J}_{\widehat{\gg}}$ are images by $\psi$ of members of the above  minimal set of generators of  $\mathcal{J}_\gg$.
Therefore, one gets $\widetilde{\mathcal{J}_\gg} : x_1x_3 =
\mathcal{J}_{\widehat{\gg}},$ hence $\mathcal{J}_{\widehat{\gg}}$ is minimal prime of
$\widetilde{\mathcal{J}_\gg}$ as a confirmation of item (c) of the theorem.

\begin{Remark}\rm
    The multiplier $x_1x_3$ in the above example is not accidental.
    Indeed, the exceptional form $x_1x_3^2y_0-x_2^3y_2$ can be expressed in terms of the members of the given set of minimal generators of $\mathcal{J}_\gg$, namely we have
    $$x_1x_3^2y_0-x_2^3y_2= x_2^3y_0- x_1^2x_3y_3 -x_1x_3(x_3y_2-x_1y_3).$$
    Since this expression is an elementary operation on the members of the given set of minimal generators of $\mathcal{J}_\gg$, this example answers negatively the question as to whether there is in fact a minimal set of generators of  $\mathcal{J}_\gg$ mapping via $\psi$ $\underline{onto}$ a minimal set of generators of $\mathcal{J}_{\widehat{\gg}}$.
\end{Remark}

\subsection{Interaction with de Jonqui\`eres maps}

The so-called plane de Jonqui\`eres maps are at the heart of
classical plane Cremona map theory and it was generalized to higher
dimensional space $\pp^n$ with $n \geq 3$. The result of this
section enhances the role of the Newton dual for these maps. We
briefly recall a few preliminaries on these maps (for more details
see \cite{PS2015}). For each given point of $\pp^n$ there is a
subgroup of the entire Cremona group of dimension $n$ consisting of
such maps. Fixing the point $o = (0:\ldots:0:1)$, these maps will
form the de Jonqui\`eres subgroup $J_o(1; \pp^n) \subset Cr(n)$ of
type $1$ with center $o$, where $Cr(n)$ is the Cremona group of
$\pp^n$. Let $H := \{x_n = 0\}$, by \cite[Proposition 1.2]{PS2015}, a
Cremona map belongs to $J_o(1; \pp^n)$ can be defined by
$$(qg_0: \ldots : qg_{n-1} : f),$$ where $(g_0 : \ldots : g_{n-1})$
defines a Cremona map of $H \simeq \pp^{n-1}$ and $q, f \in k[\xx]$
are relatively prime $x_n$-monoids one of which at least has
positive $x_n$-degree.
The Cremona map defined by $(g_0 : \ldots : g_{n-1})$ is the {\em underlying} or {\em support map} of the de Jonqui\`eres map.

Perhaps the simplest of the de Jonqui\`eres map in any dimension is the Magnus involution, where one can take $q=x_n$ and $f=x_0\cdots x_{n-1}$, while the support map is of the same kind in one dimension less.
The next result shows that de Jonqui\`eres
maps are preserved under the Newton dual transform.

\begin{Proposition}\label{PropDeJonq}
If  $\mathfrak{F} \in
J_o(1; \pp^n)$ then $\widehat{\mathfrak{F}} \in
J_o(1; \pp^n)$.
\end{Proposition}
\demo By \cite[Proposition 1.2]{PS2015}, an element of $J_o(1; \pp^n)$  can be characterized as having  a representative of the form
$$\{qg_0, \ldots  qg_{n-1},  f\},$$
where $(g_0 : \ldots : g_{n-1})$ defines a Cremona map of $H \simeq \pp^{n-1}$ and $q, f \in k[\xx]$ are relatively prime $x_n$-monoids one of which at least has positive $x_n$-degree.
Write $\ff$ for such a representative of $\mathfrak{F}$.

Let ${\boldsymbol\alpha} := (\alpha_0, \ldots, \alpha_n)^t$,
${\boldsymbol\beta} := (\beta_0, \ldots, \beta_n)^t$,
${\boldsymbol\gamma} := (\gamma_0, \ldots, \gamma_n)^t$ and
${\boldsymbol\delta} := (\delta_0, \ldots, \delta_n)^t$ stand for the
directrix vectors of $N(\ff)$, $N(\gg)$ , $N(q)$ and $N(f)$, respectively. By
Theorem~\ref{TheoChapeuN},
$$
\begin{array}{ccl}
  \widehat{\ff} & = & (x_0^{d-\alpha_0} \ldots \: x_n^{d-\alpha_n})^{-1}\ff(\widehat{\xx}) \\
                & = & (x_0^{d-\alpha_0} \ldots \: x_n^{d-\alpha_n})^{-1} \: \{q(\widehat{\xx}) \gg (\widehat{\xx}), f(\widehat{\xx})
                \} \\
                & = & (x_0^{d-\alpha_0} \ldots x_n^{d-\alpha_n})^{-1} \: \{ (x_0^{r-\gamma_0} \ldots x_n^{r-\gamma_n}) \: (x_0^{s-\beta_0} \ldots x_n^{s-\beta_n}) \: \widehat{q} \: \widehat{\gg}, (x_0^{d-\delta_0} \ldots x_n^{d-\delta_n}) \: \widehat{f}
                \} \\
                & = & \{ (x_0^{\alpha_0 - \gamma_0 - \beta_0} \ldots \: x_n^{\alpha_n - \gamma_n - \beta_n}) \: \widehat{q} \: \widehat{\gg}, (x_0^{\alpha_0-\delta_0} \ldots x_n^{\alpha_n-\delta_n}) \: \widehat{f}
                \},
\end{array}
$$
where $d := deg(\ff) = deg(f)$, $s :=
deg(\gg)$ and $r := deg(q)$. By Corollary~\ref{corobira},
$\widehat{\gg}$ is a Cremona map of $H \simeq \pp^{n-1}$. So, it suffices to prove that
$(x_0^{\alpha_0 - \gamma_0 - \beta_0} \ldots \: x_n^{\alpha_n -
\gamma_n - \beta_n}) \: \widehat{q}$ and $(x_0^{\alpha_0-\delta_0}
\ldots x_n^{\alpha_n-\delta_n}) \: \widehat{f}$ are relatively prime
$x_n$-monoids one of which at least has positive $x_n$-degree. Since
$f$ and $q$ are relatively prime $x_n$-monoids one of which at least has positive $x_n$-degree, then so are $\widehat{f}$ and
$\widehat{q}$. Note that $\beta_n = 0$, $\alpha_n = 1$ and $\alpha_i = max\{\gamma_i + \beta_i,\delta_i\}$, for every $i = 0,\ldots,n$. Therefore $(x_0^{\alpha_0 - \gamma_0 - \beta_0} \ldots \: x_n^{\alpha_n - \gamma_n - \beta_n}) \: \widehat{q}$ and $(x_0^{\alpha_0-\delta_0} \ldots \: x_n^{\alpha_n-\delta_n}) \: \widehat{f}$ are relatively prime $x_n$-monoids one of which at least has positive $x_n$-degree. \qed

\medskip

We now come to the question as which de Jonqui\`eres maps have the
property that its Newton dual commute with its inverse map. Recall that this property is equivalent to having the map commute with the Magnus involution $ \mathfrak{M}$ (Corollary~\ref{compMagnus} (ii)).
Thus, we will phrase a partial answer to the question in this language.

\begin{Proposition}\label{JonqComM}
Let $\mathfrak{F}\in J_0(1;\pp^n)$ be defined by the forms  $qg_0,\ldots, qg_{n-1}, f$, where the underlying Cremona map $(g_0 : \ldots : g_{n-1})$  of $H := \{x_n = 0\} \simeq \pp^{n-1}$ is monomial and $q, f \in k[\xx]$ are relatively prime forms.  The following conditions are equivalent:
\begin{itemize}
\item[(i)] $\mathfrak{M} \circ \mathfrak{F} = \mathfrak{F} \circ \mathfrak{M};$
\item[(ii)] $\widehat{q'} = f'$,
\end{itemize}
where $q := q'  M_q$  and $f := f' M_f$ with $M_q$ and $M_f$  monomials of highest degrees  dividing $q$ and $f$, respectively.
\end{Proposition}
\demo By Corollary~\ref{compMagnus}  (i), $\mathfrak{F} \circ \mathfrak{M} = \widehat{\mathfrak{F}}$. In addition, by the proof of Proposition~\ref{PropDeJonq}, the Newton dual $\widehat{\mathfrak{F}}$ is defined by the forms
$$\{T_1 \: \widehat{q} \: \widehat{\gg}, T_2 \widehat{f}\},$$
where $T_1$ and $T_2$ are suitable monomials and this representative has no fixed part. On the other hand, $\mathfrak{M} \circ \mathfrak{F}$ admits the representative $$\{f g_1 \cdots g_{n-1}, f g_0 g_2 \cdots g_{n-1},\ldots, f g_1 \cdots g_{n-2}, q g_0 \cdots g_{n-1}\}.$$
Let $T$ be the fixed part of representative above, i.e., the $\gcd$ of its terms. Since $\{g_0,\ldots,g_{n-1}\}$ are monomials and gcd$(q,f) = 1$, $T$ is a monomial.

\medskip

$(i) \Rightarrow (ii)$.  Drawing on the format of the preceding representatives the assumed equality $\mathfrak{M} \circ \mathfrak{F} = \mathfrak{F} \circ \mathfrak{M}$ implies that $T T_2 \widehat{f} = q g_0 \cdots g_{n-1}$ which, by Corollary~\ref{ComNewMult}, implies that $\widehat{\widehat{f}} = \widehat{q}$ and hence $f' = \widehat{q'}$, by the usual behavior while taking the Newton dual of a form.

\medskip

$(ii) \Rightarrow (i)$. We argue that the above respective representatives of  $\mathfrak{M} \circ \mathfrak{F}$ and $\mathfrak{F} \circ \mathfrak{M}$ are proportional, hence showing they define the same map.

In other words, we claim that the following matrix has rank $1$
$$A:=\left(
                \begin{array}{ccccc}
                  T_1 \: \widehat{q} \: h_0  & T_1 \: \widehat{q} \: h_1 &\ldots & T_1 \: \widehat{q} \: h_{n-1} & T_2 \: \widehat{f}\\
                   f g_1 \cdots g_{n-1} & f g_0 g_2 \cdots g_{n-1} & \ldots & f g_1 \cdots g_{n-2} & q g_0 \cdots g_{n-1}\\
                 \end{array}
     \right)
,$$
where $\widehat{\gg} = \{h_0,\ldots,h_{n-1}\}$.
Let $A_{i,j}$ denote the typical $2\times 2$ submatrix of $A$ with the $i$th and $j$th columns ($0\leq i < j\leq n$). For $j \leq n-1$, one has
\begin{eqnarray}\nonumber
\det(A_{i,j})&=&(T_1 \: \widehat{q} \: h_i)\left( f \: \frac{g_0 \cdots g_{n-1}}{g_j}\right) - (T_1 \: \widehat{q} \: h_j)\left( f \: \frac{g_0 \cdots g_{n-1}}{g_i}\right)\\ \nonumber
&=& T_1 \: \widehat{q} \: f \: \frac{g_0 \cdots g_{n-1}}{g_ig_j}(h_ig_i - h_jg_j)\\ \nonumber
& =& T_1 \: \widehat{q} \: f \: \frac{g_0 \cdots g_{n-1}}{g_ig_j}\left(x_0^{\beta_0} \ldots \: x_n^{\beta_n} - x_0^{\beta_0} \ldots \: x_n^{\beta_n}\right)\\ \nonumber
&=& 0,
\end{eqnarray}
where ${\boldsymbol\beta} := (\beta_0, \ldots, \beta_n)^t$ is the directrix vector of $N(\gg)$.
For $j=n$, one finds
\begin{eqnarray}\nonumber \det(A_{i,n})&=&(T_1 \: \widehat{q} \: h_i)(q g_0 \cdots g_{n-1}) - (T_2 \: \widehat{f})\left( f \: \frac{g_0 \cdots g_{n-1}}{g_i}\right)\\ \nonumber
&=& \left(\frac{g_0 \cdots g_{n-1}}{g_i}\right) (T_1 \: q \: \widehat{q} \: h_i \: g_i - T_2 \: f \: \widehat{f}).
\end{eqnarray}
By proof of Proposition~\ref{PropDeJonq}, $T_1 = x_0^{\alpha_0 - \gamma_0 - \beta_0} \ldots \: x_n^{\alpha_n - \gamma_n - \beta_n}$ and $T_2 = x_0^{\alpha_0-\delta_0} \ldots x_n^{\alpha_n-\delta_n},$ where ${\boldsymbol\alpha} := (\alpha_0, \ldots, \alpha_n)^t$, ${\boldsymbol\gamma} := (\gamma_0, \ldots, \gamma_n)^t$ and ${\boldsymbol\delta} := (\delta_0, \ldots, \delta_n)^t$ are the directrix vectors of $N(\ff)$, $N(q)$ and $N(f)$, respectively.
Therefore, we get
$$
\begin{array}{ccl}
  T_1 \: q \: \widehat{q} \: h_i \: g_i & = & x_0^{\alpha_0 - \gamma_0} \ldots \: x_n^{\alpha_n - \gamma_n} q \: \widehat{q}\\
                                        & = & x_0^{\alpha_0 - \gamma_0} \ldots \: x_n^{\alpha_n - \gamma_n} M_q \: q' \: \widehat{M_q q'} \\
                                        & = & x_0^{\alpha_0 - \gamma_0} \ldots \: x_n^{\alpha_n - \gamma_n} M_q \: q' \: \widehat{q'} \\
                                        & = & x_0^{\alpha_0 - \gamma_0} \ldots \: x_n^{\alpha_n - \gamma_n} M_q \: q' \: f'
\end{array}
$$
and
$$
\begin{array}{ccl}
  T_2 \: f \: \widehat{f} & = & x_0^{\alpha_0-\delta_0} \ldots x_n^{\alpha_n-\delta_n} \: f \: \widehat{f} \\
                          & = & x_0^{\alpha_0-\delta_0} \ldots x_n^{\alpha_n-\delta_n} \: M_f f' \: \widehat{M_f f'} \\
                          & = & x_0^{\alpha_0-\delta_0} \ldots x_n^{\alpha_n-\delta_n} \: M_f f' \: \widehat{f'} \\
                          & = & x_0^{\alpha_0-\delta_0} \ldots x_n^{\alpha_n-\delta_n} \: M_f q' f'.\\
\end{array}
$$
Write $M_q: = x_0^{b_0} \ldots x_n^{b_n}$ and $M_f: = x_0^{c_0} \ldots x_n^{c_n}$. Since $\widehat{q'} = f'$, $\widehat{\widehat{q'}} = q'$ and $\widehat{\widehat{f'}} = f'$, then $q'$ and $f'$ have the same directrix vector ${\boldsymbol\eta} := (\eta_0, \ldots, \eta_n)^t$. In addition, $q = M_q q'$  and $f = M_f f'$, hence $\boldsymbol\eta = \boldsymbol\gamma - (b_0,\ldots,b_n)$ and $\boldsymbol\eta = \boldsymbol\delta - (c_0,\ldots,c_n)$.
Therefore,  $\det(A_{i,n})$ vanishes. \qed


\end{document}